\documentclass{article}
\usepackage{graphicx}
\usepackage{amsmath}
\usepackage{amsfonts}
\usepackage{amssymb}
\usepackage[all]{xy}
\usepackage{natbib}
\usepackage[latin1]{inputenc}
\usepackage{float}
\newtheorem{theorem}{Theorem}

\newtheorem{Definition}{Definition}
\newtheorem{example}{Example}
\newtheorem{result}{Result}
\newtheorem{proposition}{Proposition}

\newcommand{\ind}[1]{1\hskip -2.5 pt\hbox{I}_{\{#1\}}}

\begin{document}

\title{\textbf{Properties of Doubly Stochastic Poisson Process with affine intensity}}
%\author{\textit{Alan De Genaro Dario^{ $\S$ $\ddag$}  \\ {\footnotesize $\ddag$ Institute of Mathematics and Statistics  - IME$\backslash$USP}and Adilson Simonis $\ddag$}\\ {\footnotesize  $\S$ Securities, Commodities and Futures Exchange - BM\&FBOVESPA}}

\date{}
\author{\textit{Alan De Genaro Dario $^{\S \ddag}$\footnote{The results of this paper are part of the first author's Ph.D Thesis completed under supervision of the second author. Alan De Genaro Dario would like to thank Marco Avellaneda, Jorge Zubelli, Cristiano Fernandes and Julio Stern, members of his Ph.D Thesis committee, and Yuri Suhov for comments and suggestions - \emph{e-mail:} adario@bvmf.com.br} \hspace{.2cm} and \hspace{.2cm} Adilson Simonis $^{\ddag}$\footnote{\emph{e-mail:} asimonis@ime.usp.br} }\\ {\footnotesize $\S$ Securities, Commodities and Futures Exchange - BM\&FBOVESPA} \\ {\footnotesize $\ddag$ Institute of Mathematics and Statistics  - IME$\backslash$USP}}
\maketitle

\begin{abstract} \noindent {This paper discusses properties of a Doubly Stochastic Poisson Process
(DSPP) where the intensity process belongs to a class of affine
diffusions. For any intensity process from this class we derive an
analytical expression for probability distribution functions of
the corresponding DSPP. A specification of our results is provided
in a particular case where the intensity is given by
one-dimensional Feller process and its parameters are estimated by
Kalman filtering for high frequency transaction data.}

% Additionally the stationarity of
%the point process is addressed. To illustrate our results in this
%paper some common affine diffusion used to modelling interest rate
 %term structure are adopted to drive the Doubly Stochastic Poisson
%Process

\end{abstract}

\section{Introduction}

The Doubly Stochastic Poisson Processes (DSPPs) were introduced by
\citet{Cox1} by allowing the intensity of the Poisson process to
be described by a positive random variable and not just a
deterministic function. The aim of such a generalization was to
allow the dynamic of a process that is exogenous to the model to
influence transitions in the point process that we are concerned
with. Point processes have applications in several areas of
applications, among which we mention: biostatistics, finance and
reliability theory. In biostatistics they form a theoretical
framework for studying recurrent events, as was done in
\citet{gali}, in studies of the size of tumors in rats over a
period of time. In finance, \citet{lando} was the pioneer in using
point processes for describing occurrences of credit events. In
reliability theory, \citet{dalal} work has developed criteria for
determining an optimum stopping time for testing and validating a
software. \\

%For the time being, and without the necessary formalism that will
%be adopted in the next section, the $N_t$ point process can be
%interpreted as the representation of the act of counting a
%phenomenon over time.

% In the following sections we will make a brief review of
%the principal work done on this subject, as well as describing the
%problem and aims of this article.

It was the seminal work by \citet{Cox1} which introduced the DSPPs
(also known as the Cox Process). The main work in this field is
\citet{Grandell} where the main properties of the DSPPs have been
presented in terms of standard construction of Probability Theory.
Alternatively, in the books of \citet{Bremaud} and \citet{Jones}
the presentation is based in line with the concepts and properties
of Martingales. However, all these sources focus on deriving
general properties of DSPP, without exploring  the functional form
of the intensity of the process. Consequently, they attracted a
reduced number of applications.\\

The study of the DSPPs took a new turn once the functional form
for the intensity of the process has been specified. As a result
it became possible to obtain analytical expressions for
probability density functions for different types of processes. In
this context we can quote  may be cited \citet{Bouzas1} who made
use of truncated normal distribution to describe the intensity of
the process and \citet{Bouzas2} who generalized the form of
intensity to include the case of a harmonic oscillator. The
contribution of these works was that closed analytical expressions
for density functions of the Cox process have been obtained, as
well as their moments. However, in both cases the authors used
constructions in which the restriction on non-negativity for the
intensity measure was not maintained. In order to get a round of
this limitation, the authors defined a region in the parameter
space where the probability of occurrence of negative values for
the intensity is reduced. With the intent to guarantee the
preservation of the non-negativity condition, \citet{BasuDassios}
and \citet{Kozachenko} suggested the adoption of a lognormal model
for the process intensity. A formulation that incorporates the
intensity into a dynamic formulation is developed in
\citet{DassiosJang} who use the functional form of a process of
the Shot-Noise type which, despite guaranteeing the non-negativity
of intensity, is not in fact a diffusion process. On the other
hand \citet{wei} assumed that the intensity is governed by a
one-dimensional Feller process and
obtained a form of probability density function for the corresponding DSPP.\\

Feller processes were introduced and established in the financial
literature after \citet{CIR}. One of the properties is that the
Feller process lives in $\mathbb{R}_+$ which guarantees that the
non-negativity condition is fulfilled. In the present work we
assume that the intensity is controlled by a related diffusion
process, as formalized by \citet{DK}, which incorporates Feller
processes in one or $d$ dimensions. In this way, the models from
\citet{wei}, \citet{BasuDassios} and \citet{Kozachenko}
can be seen as specific examples of the model proposed in the present work.\\

The use of point processes in finance, especially of DSPPs,
progressed considerably at the end of the 1990s with the
development of models of managing and pricing the credit risk. In
particular, \citet{DS} and \citet{DFS} formalized the construction
of the probability density of the first jump in the process, in a
related context. More precisely, the time to the first jump
represents, in the context of credit risk, the time until the
bankruptcy (default) of a company (and/or a country). Here, once
the absorbing state had been reached, it was unnecessary to study
the further dynamic of the DSPP. \\

Recently a new area of applications of point processes in finance
has emerged, along with the use of these models for describing the
arrival process of bid and ask orders in an electronic trading
environment. In these models, the arrival process of orders
changes over time;  the idea is to characterize the dynamic of
this process and obtain expressions that may be treated
analytically, describing the probability that an order had been
sent in line with a market configuration and was executed before
the price was
altered. See \citet{RamaCont}.\\

Perhaps the most well known point process is the homogenous
Poisson process. For this process the arrival rate is constant. A
homogenous Poisson process can therefore be described by a single
value $\lambda_t =\lambda$. However, in many applications the
assumption of a constant arrival rate is not realistic. Indeed, in
financial data we tend to observe bursts of trading activity
followed by lulls. This feature becomes apparent when looking at
the series of actual orders arrivals. Figure 1 presents the
average number of sell orders for BRL/USD FX futures submitted to
the Brazilian Exchange, BVMF. The plot indicates a non-constant
behavior for the average number of submitted orders; the
homogenous Poisson model is clearly not suitable for such data.

%\vspace{-1.2cm}
\begin{figure}[!h]
\centering
\includegraphics[scale=0.4]{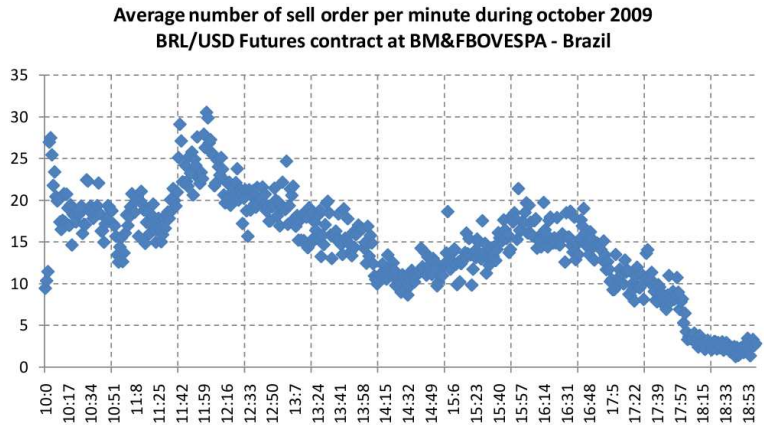}\vspace{-1.5cm} %{ofertas4.pdf}
\caption{Sell orders for FX futures contract during October:
(ticker DOL FUT - expiry NOV09)}
\end{figure}

In fact, the number of bids and asks that getting into the order
book may depend upon a number of factors exogenous to the model.
For example, it may be the level of investors risk aversion on a
given day, or intra-day seasonality, or a disclosure of a piece of
information, or a new technology producing an impact on a certain
sector of the economy. For this reason, in our opinion, the
construction of a model that can be treated analytically and
incorporates endogenously a stochastic behavior of the intensity
of a point process represents a contribution to the literature. In
light of the above topics, the main aim of this article is to
present new results on DSPPs when the intensity
belongs to a family of so-called affine diffusions with potencial application in high frequency trading.\\

In this situation seems reasonable to focus on a particular class
of models and attempt to use the approach based on point processes
in order to study their dynamic over time. For this purpose we
selected the Affine Term Structure (ATS) models, and our goal is
to obtain the form Theorem for the Probability Density Functions
for the Cox Process when the intensity belongs to a family of
affine diffusion. Additionally, in one particular case when the
intensity is governed by an one-dimensional Feller diffusion we
obtain more detailed results, such as its moments and the
convergence to its stationary distribution. Finally we propose an
estimation procedure for point processes with stochastic affine
intensity based on Kalman Filter conjugated with quasi-maximum
likelihood estimator. Thus, the estimation procedure is applied to
high frequency transactional data from FX futures contracts traded
in Brazil.

%======================================

\section{The basic structure}
Consider the filtered probability space $\left( \Omega
,\mathcal{G},\mathbb{G}, \mathbb{P}\right)$, where
$\mathbb{G}=(\mathcal{G}_{t})_{t\geq 0}$ is a filtration having
the sets with null measure and right continuous.

\begin{Definition}
Let $\tau$ be a non-negative random variable
$\tau:\Omega\rightarrow\mathbb{R}^+\cup\infty$, which is defined
on a probability space $\left( \Omega ,\mathcal{G} , \mathbb{G},
\mathbb{P}\right)$, such that $\{\tau\leq t \} \in \mathcal{G}_{t}
\hspace{.3cm} \forall t \geq 0$.
\end{Definition}

\begin{Definition} Given a no-decreasing sequence of stopping time $\{\tau_i, i
\in\mathbb{N}\}$, the \textbf{Point Processes} $N_t$ is a
(\textit{càdlàg}) processes defined as:
\end{Definition}
\begin{equation}
    N_t:= \displaystyle\sum_{i}\ind{\tau_i\leq t}
\end{equation}

\begin{Definition} Define the \textbf{intensity processes} as:
\vspace{-.3cm}
\begin{align}\label{inten}
    \lambda_t &=\rho_0+\sum_{i=1}^d \rho_{i,1} X_{i,t}\\
              &=\rho_0+\rho_1 \cdot X_t
\end{align}
With $\rho_0 \in \mathbb{R}$ and $\rho_1 \in \mathbb{R}^{d}$.
Where the state variable $X_t$ will be defined below.\\
\end{Definition}

\begin{Definition}
Let $\left( \Omega ,\mathcal{G},\mathbb{G}, \mathbb{P}\right)$ be
a filtered probability space and let $\mathbb{F}$ be a
sub-filtration of $\mathbb{G}$. We called \textbf{cumulate
intensity} (or \textbf{Hazard Process}) an $\mathcal{F}$-adapted,
right-continuous, increasing stochastic process
$\Lambda_{0,t}=\int_0^t\lambda_u du$ for $t\geq 0$ with
$\mathbb{P}$-a.s. $\Lambda_0=0$ and $\Lambda_\infty=\infty$.\\
\end{Definition}

The hazard process plays an important role in the martingale
approach to either credit risk or asset pricing because it is the
compensator of the associated doubly stochastic Poisson process,
see for example Bielecki and Rutkowski (2002). For our purposes,
the hazard process will be important to derive the probability
density function for the Point Process, $(N_t)_{t \geq 0}$. To
construct the stopping time with the given intensity, define as in
Lando (1998):

\begin{equation}
    \tau_i:=\inf \Bigl\{t: \Lambda^i_{0,t} \geq E_{1,i}\Bigr\}
\end{equation}
where $E_{1,i}$ is a unit exponential random variable independent
of $X_t$ $\forall$ \emph{i}.\\

\begin{Definition} Let $X_t:\mathcal{D}\times\mathbb{R}
\rightarrow\mathbb{R}^d$, be the \textbf{state variable} as
solution for the stochastic differential equation SDE:
\end{Definition}

\begin{equation}\label{SDE}
dX_{t} =\mu(X,t)dt+\sigma(X,t)dW_{t}
\end{equation}

\noindent where: $\mu(X,t) :\mathcal{D} \times \mathbb{R_+}
\rightarrow \mathbb{R}^d$, is the drift; $\sigma(X,t) :
\mathcal{D} \times \mathbb{R}_+ \rightarrow \mathbb{R}^{d\times
d}$, the diffusion coefficient and, $W_{t} $ the standard Brownian
Motion defined in $\left( \Omega ,\mathcal{G},\mathbb{G},
\mathbb{P}\right)$. Where $\mathcal{D}\subset\mathbb{R}^d$ will be
properly defined below.
\\

\begin{Definition}[\citet{DK}]\label{SDE2}
Define $X_t$ as an \textbf{affine processes}, if the following
condition are simultaneously satisfied for (\ref{SDE}):
\end{Definition}

\begin{enumerate}
    \item Drift
\begin{equation}\label{mu}
    \mu(X,t)=\mathcal{K}(\Theta-X_t)
\end{equation}

For $\Theta \in \mathbb{R}^d$ and $\mathcal{K} \in
\mathbb{R}^{d\times d}$

    \item Covariance matrix

\begin{equation}\label{sig}
    \sigma(X,t) = \Sigma\sqrt{\sigma}
\end{equation}
%For $H_{0ij} \in \mathbb{R}^d$ e $H_{1ij} \in \mathbb{R}^{d}$
\end{enumerate}

Where $\Sigma$ is $d \times d$ matrix not necessarily symmetric,
and $\sigma$ is a diagonal matrix with the \emph{i}-th elements
from the diagonal given by:

\begin{equation}\label{afim_vol}
\sigma_{i}= a_{i}+b_{i} X_t, \hspace{.5cm} i=1,\ldots,d,
\end{equation}

with \emph{a} $ \in \mathbb{R}$ e $b \in \mathbb{R}^d$.\\

Adopting the construction of \citet{DK} and \citet{dai} the state
variable in its affine form is:

\begin{equation}\label{afim}
    \mbox{d}X_t = \mathcal{K}(\Theta-X_t)dt+\Sigma\sqrt{\sigma}\mbox{d}W_t
\end{equation}

where $W_t$ is the \emph{d-dimensional} independent standard
Brownian Motion in $\left( \Omega ,\mathcal{G},\mathbb{G},
\mathbb{P}\right)$, $\mathcal{K}$ and $\Sigma$ are $d \times d$
matrices which may be nondiagonal and
asymmetric.\\

According to \citet{DK} the coefficient vectors $b_{1}, \ldots,
b_{d}$ in (\ref{afim_vol}) generate \emph{stochastic volatility}
unless they are all zero, in which case (\ref{afim}) defines a
Gauss-Markov process. If the volatilities (all or some) are
stochastic, then two questions face up: how to constrain the model
in order to avoid negative volatilities? Under those constraints,
which is the set $\mathcal{D}\subset\mathbb{R}^d$ where the
process $(X_t)_{t \geq 0}$ can take values?
\newline
The open domain $\mathcal{D}$ implied by nonnegative volatilities
may be defined as:

\begin{equation}\label{adm3}
    \mathcal{D} = \bigg\{ x \in \mathbb{R}^d : a_{i}+b_{i}\cdot x
    \geq 0, \hspace{.5cm} i \in \{1,\ldots,d\}\bigg\}
\end{equation}

The following condition of \citet{DK} is sufficient for the
positiveness of $\sigma_i$:

\textbf{Condition A (\citet{DK})}: \textit{For all i}:
\begin{enumerate}
    \item $\forall x \hspace{.3cm} \mbox{such that} \hspace{.3cm} \sigma_i(x)=0, \hspace{.3cm} b_i^\top(\mathcal{K}(\Theta-X_t)>\frac{1}{2} b_i^\top\Sigma\Sigma^\top b_i$
    \item $\forall j,\hspace{.3cm} \mbox{if} \hspace{.3cm} (b_i^\top \Sigma)_j \neq 0$,
    then $\sigma_i = k \sigma_j$ for $k>0$
\end{enumerate}

Both parts of Condition A are designed to ensure strictly positive
volatility, and they are both effectively necessary for this
purpose.\\

A second set of constrains must be imposed to guarantee that the
non-negativity condition for $(\lambda_t)_{t \geq 0}$ is
fulfilled. The open domain $\mathcal{D}$ implied by nonnegative
intensities may be defined as:

\begin{equation}\label{adm4}
    \mathcal{D} = \bigg\{ x \in \mathbb{R}^d : \lambda_t(x) \geq 0\bigg\}
\end{equation}

In this sense, the strongest form of guarantee that (\ref{adm3})
and (\ref{adm4}) are simultaneously met it is to impose
$\mathcal{D}=\mathbb{R}^d_+$. However, it is a very strong
constraint on the state space resulting in a large number of
models excluded.\\

Using the Canonical Representation developed by \citet{DS} it will
be possible to overcome the non-negativity for $(\lambda_t)_{t
\geq 0}$. In fact, according to the Canonical Representation,
$\mathbb{A}_m(d)$, where $m=rank(\sigma)$, the state vector $X_t$
will be split into two subvectors: the first factors vector
$X_{t}^B \in \mathbb{R}^m$, with $0 \leq m \leq d$ and the other
group of factors are stacked in $X_{t}^D
\in \mathbb{R}^{d-m}$.\\

This way it is also possible to overcome existence and uniqueness
problems for a more general class of affine diffusion than the one
that satisfies (\ref{afim}). In fact the process $X_{t}^B$ exists
and is unique, moreover it is autonomous with respect to
$X_{t}^D$. The volatility of $X_{t}^D$, conditionally on
$X_{t}^B$, is given, so there is no uniqueness problem concerning
$X_{t}^D$. Finally, from the Canonical Representation it is known
that the process $X^D$, when $X^D\neq 0$, has Gaussian
Distribution, independent from $X^B$, with mean $\mu^D$ and standard deviation $\sigma^D$.\\

From the Canonical Representation, we just need to impose an extra
condition to ensure strictly positive of $(\lambda_t)_{t \geq
0}$.\\

\textbf{Condition B}: if we set
\begin{equation}\label{cond_5}
    \mu^D > 0 \hspace{.1cm},
\end{equation}

it is possible to state that:
\begin{proposition}\label{positivo} If $\gamma^{*}\mu^D > \sigma^D$, then the process $X^D$
is nonnegative almost sure, i.e. $X^D \in \mathcal{D}^+$.
\end{proposition}

\noindent {\bf Proof of Proposition \ref{positivo}}\\

According to the Canonical Representation, the process $X^D$ is
Gaussian  with mean $\mu^D$ and standard deviation $\sigma^D$.
Thus it is possible to compute the probability $X^D \in
\mathcal{D}^+$:

\begin{equation}\label{positivo_2}
    \mathbb{P}(X^D > 0)= \Phi \left(\frac{\mu^D}{\sigma^D}\right)
\end{equation}

where $\Phi(\cdot)$ is the cumulative normal distribution function.\\

Because condition (\ref{cond_5}) the problem described in
(\ref{positivo_2}) is now written as:

\begin{equation}\label{positivo_3}
 \gamma^{*} =  \displaystyle \inf_\gamma(\gamma:\Phi(\gamma)=1)
\end{equation}
where: $\gamma= \frac{\mu^D}{\sigma^D}$\\

Then by requiring that $\gamma^{*}\mu^D > \sigma^D$ the
non-negativity condition for $(\lambda_t)_{t \geq 0}$ is satisfied
almost sure.

\begin{flushright}
$\blacksquare$
\end{flushright}

We now turn to examples and illustrations of Affine framework.
\begin{itemize}
    \item \begin{example}\hspace{-.12cm}\emph{:} Ornstein-Ulhenbeck (Vasicek)
$$dX_t=\kappa(\theta-X_t)dt+\sigma dW_t$$
\end{example}
    \item \begin{example}\hspace{-.12cm}\emph{:} Feller (or Square Root, Cox-Ingersoll-Ross)
$$dX_t=\kappa(\theta-X_t)dt+\sigma\sqrt{X_t} dW_t$$
\end{example}
    \item \begin{example}\hspace{-.12cm}\emph{:} Geometric Brownian Motion
$$\frac{dS_t}{S_t}= \mu dt+\sigma dW_t$$
$$X_t:=ln S_t \hspace{1cm} dX_t=\left(\mu-\frac{1}{2}\sigma^2\right)dt+\sigma
dW_t$$
\end{example}
    \item \begin{example}\hspace{-.12cm}\emph{:} Heston (Stochastic Volatility)
    $$X_t:=ln S_t \hspace{1cm} dX_t=\left(\mu-\frac{1}{2}\sigma^2_t\right)dt+\sigma_t
    dW^1_t$$ \\
$$d\sigma_t=\kappa(\theta-\sigma_t)dt+\nu\sqrt{\sigma_t} dW^2_t$$\\
$$<dW^1_t,dW^2_t>=\rho$$
\end{example}
    \item \begin{example}\hspace{-.12cm}\emph{:} Multivariate CIR
$$
dX_t = \left(\left[\begin{array}{c}
  \kappa_1\theta_1 \\
  \kappa_2\theta_2 \\
\end{array}\right]+\left[
\begin{array}{cc}
  -\kappa_1 & 0 \\
  0 & -\kappa_2 \\
\end{array}
\right]
 X_t\right)dt + \left[
\begin{array}{cc}
  \sigma_1 & 0 \\
  0 & \sigma_2 \\
\end{array}
\right] \left(\sqrt{\left[
\begin{array}{cc}
  1 & 0 \\
  0 & 1 \\
\end{array}
\right]X_t}\right)dW_t
$$
\end{example}
\end{itemize}

%A filtração com as informações suficiente para descrever a
%dinâmica do processo de intensidade $\lambda_t$ é a
%sub-$\sigma$-algebra de $\mathcal{G}_t$, definida por
%$\mathcal{F}_t :=\sigma\{W_s:0\leq s \leq t\}$.

In order to have the no-negative random variable $\tau_i$ as a
stopping time a technical conditions must be imposed over the
filtration. So we shall first examine a trivial case where the
condition is not respected.

\begin{proposition}\label{filta} The hitting time $\tau$ \textbf{is not} a $\mathbb{F}$-stopping time with respect to
$\mathcal{F}_t :=\sigma\{W_s:0\leq s \leq t\}$
\end{proposition}

\noindent {\bf Proof of proposition \ref{filta}} \\
Assume, by absurd, that $\tau$ is $\mathbb{F}$-stopping time for
$\mathcal{F}_t :=\sigma\{W_s:0\leq s \leq t\}$, so using the
Martingale Representation Theorem there exist a compensated
process (a $\mathbb{F}$-Martingale) $M_t$ which may be represented
as a stochastic integral with respect to the Brownian motion $W$.
Therefore we arrive a contradiction because $M_t$ must jump in
$\tau$.
\begin{flushright}
$\blacksquare$
\end{flushright}

Thus the filtration $\mathbb{F}$ should be enlarge. There are
several ways to expand it but we do enlarged $\mathbb{F}$ just to
get $\tau$ as stopping time. Then the proper filtration
$\mathbb{G}=(\mathcal{G}_{t})_{t\geq 0}$ will be constructed as:

\begin{equation*}
\mathcal{G}_t = \mathcal{F}_t \vee\mathcal{H}_t
\end{equation*}
where\footnote{$A \vee B$:= maximum between A and B}
\begin{equation*}
\mathcal{H}_t =\sigma\{N_s :0\leq s \leq t\}
\end{equation*}

\begin{Definition}\label{PPDE}

Let $\left( \Omega ,\mathcal{G},\mathbb{G}, \mathbb{P}\right)$ be
a filtered probability space and $N_t$ a stochastic processes
defined on it. Let $\mathbb{F}$ be a sub-filtration of
$\mathbb{G}$, then $N_t$ is called a Doubly Stochastic Poisson
Process with relation to $\Lambda_t$, if $\lambda_t$ is
$\mathbb{F}$-measurable  and for all $0\leq t \leq T$ and
$k=0,1,\ldots$ the next condition be satisfied:
\begin{equation}\label{pdf}
P(N_T-N_t=k|\mathcal{H}_t\vee \mathcal{F}_\infty)=
\frac{1}{k!}\left(\Lambda_{t,T}\right)^k
\exp\left(-\Lambda_{t,T}\right)
\end{equation}
%$\mbox{para} \hspace{.3cm}  k=0,1,2,...$
\\
where $\mathcal{F}_\infty=\sigma(\mathcal{F}_u: u \in
\mathbb{R}_+)$
\end{Definition}
As a particular case:
\begin{equation}\label{pdf3}
P(N_T-N_t=k|\mathcal{H}_t\vee \mathcal{F}_\infty)= P(N_T-N_t=k|
\mathcal{F}_\infty)
\end{equation}

Therefore, conditioned to $\sigma$-field $\mathcal{F}_\infty$ the
increments of $N_T-N_t $  are independent of the $\sigma$-field
$\mathcal{H}_t$.\\

Thus, the probability of no occurrences within the interval
$[t,T]$ for the processes $N_t$ with intensity $\lambda_t>0$ is
given by:

\begin{equation}
P(N_T-N_t=0)=\mathbb{E}_{t} \left[ \exp\left(-\int_t^T\lambda_u
du\right)\right]
\end{equation}

%================================
\subsection{Laplace transform for the process $\Lambda_t$}

The Laplace transform of stochastic processes is a key ingredient
to achieve our results. However Laplace transform for integral of
stochastic processes may be obtained in closed form only for a
limited number of processes.\\

\citet{Albanese} formalize the criteria to define which processes
have a analytic form for its Laplace transform. Thus take a
diffusion $(X_t)_{t \geq 0}$ defined on $\left( \Omega,
\mathcal{G}, \mathbb{G}, \mathbb{P}\right)$ and consider the
Laplace transform defined by:

\begin{equation}
L_{t}(\mu,X,t):=\mathbb{E}_{t} \left[ f(X,s)\exp \left(-\mu
\int_t^T \phi(X,s) ds \right) \bigg |\mathcal{G}_t\right]
 \end{equation}

where $t \leq T$, $\mu \in \mathbb{C}$ and $f, \phi :
\mathbb{R} \mapsto \mathbb{R}$ two Borel functions.\\
\newline
Thus, it is possible to state that: %\vspace{.5cm}
\begin{result} [\citet{Albanese}]\label{lapla0} The class of stochastic processes with Laplace transform for its integral is given by:
\begin{equation}\label{lapla2}
dX_t=2\frac{h'(X_t)}{h(X_t)}\frac{A(X_t)^2}{R(X_t)}dt+\frac{\sqrt{2}A(X_t)}{\sqrt{R(X_t)}}dW_t
\end{equation}
With additional conditions:
\begin{enumerate}
    \item Three second order polyonymous in x: $A(x)$,$R(x)$,$Q(x)$. Such that $A(x)$ belongs to the set o
    $\{1,x,x(1-x),x^2+1\}$ e $R(X_t)\geq 0$;
    \item The function $h(x)$ is a linear combination of hypergeometric functions of the confluent \footnote{A hypergeometric function in its general form may be written as:
    $$_pF_q(\alpha_1,\ldots,\alpha_p;\gamma_1,\ldots,\gamma_q;z)$$
    For $p \leq q+1, \gamma_j \in \mathbb{C}\setminus\mathbb{Z}_+$ been represented using
    Taylor's expansion around $z=0$
    $$_pF_q(\alpha_1,\ldots,\alpha_p;\gamma_1,\ldots,\gamma_q;z)=\sum_0^\infty \frac{(\alpha_1)_n \ldots (\alpha_p)_n}{(\gamma_1)_n \ldots (\gamma_q)_n}\frac{z^n}{n!}$$}

    From type $_1F_1$ if $A(x) \in
    \{1,x\}$ and gaussian from type $_2F_1$ otherwise.
\end{enumerate}
Thus the Laplace Transform $L_{t}(\mu,X,t)$ is defined by:
\vspace{.2cm}
\begin{equation}
\phi(x)=\frac{Q(x,\mu)}{\mu R(x)}
\end{equation}
\end{result}
\vspace{.4cm}

\begin{example}\hspace{-.12cm}\emph{:   As an application from the above result we shall prove that
the Laplace transform of $\Lambda_t:= \int_0^t \lambda_u du$ when
the intensity follows a one-dimension Feller process exists.}\\
\end{example}
Thus assume that the polynomials are defined as:
\begin{equation}
A(x)= x, \hspace{.5cm} R(x) = \frac{2x}{\sigma^2}, \hspace{.5cm}
h(x)=x^{a/\sigma^2}e^{-\frac{b}{\sigma^2}x}
\end{equation}
Substituting the polynomial into (\ref{lapla2}) with a variable
change we have:
$$d\lambda_{t} =\kappa (\theta -\lambda_{t}
)dt+\sigma\sqrt{\lambda_{t} } dW_{t}$$

So in those cases where the Laplace transform $L_{t}(\mu,X,t)$
exists it is possible to apply the well know
result\footnote{Details might be found in \citet{kar}}:

\begin{result}[Feynman-Kac]\label{FK} Let $(X_t)_{t \geq 0}$ be a diffusion process with infinitesimal generator $\mathcal{A}$.
 Assume that $F \in C^{2,1}(\mathcal{D}^d\times [0,t))$ and $V
\in C^1$ is bounded. Then:

\begin{equation}\label{edp01}
F(x,t,\mu)=\mathbb{E}_{t} \left[ f(X_t)\exp \left(-\mu\int_t^T
V(x,s) ds \right) \bigg |\mathcal{G}_t\right]
 \end{equation}

it is solution of the partial differential equation (PDE):
%
%%(4.14)\tab \tab \tab \tab
\begin{equation}\label{edp1}
\left\{\begin{array}{rcl}
 \frac{\partial F}{\partial t} &=&
\mathcal{A}F(x,t,\mu)-V(x,t)F(x,t,\mu))\\
 F(x,0,\mu)&=& f(x)\hspace{.5cm} x \in \mathcal{D}
\end{array}
\right.
\end{equation}

and

\begin{equation}\label{gerador}
    \mathcal{A}F(x,t) = \frac{\partial F }{\partial t} + \frac{\partial F }{\partial x}\mu(x,t) +\frac{1}{2}\mbox{Tr}\bigg[
    \sigma(x,t)\sigma(x,t)^\top \frac{\partial ^{2} F(t,x) }{\partial X \partial X}\bigg]
\end{equation}

Furthermore \footnote{We write $Tr$ for trace of a matrix.}, if
$U(x,t,\mu) \in C^{2,1}$ solve (\ref{edp1}), then
$U(x,t,\mu)=F(x,t,\mu)$.
\end{result}

Taking the Theorem \ref{FK} it is possible to state that:

\begin{proposition}\label{sol_0} Let $(X_t)_{t \geq 0}$ be a diffusion processes satisfying the regularity condition so the solution to (\ref{edp1}) has the form:
\begin{equation}\label{sol}
    F(x,t)=e^{\alpha(t)+\beta(t)\cdot x}
\end{equation}
where the coefficients $\alpha(t)$ and $\beta(t)$ are
deterministic and satisfying\footnote{where $^{\shortmid}$ stand
for the derivative with respect to $t$.}:
\end{proposition}

\begin{align}\label{ricati}
 \beta^\shortmid(t)&=\rho_1-\mathcal{K}^\top\beta(t)-\frac{1}{2}\beta^\top b \beta \\
\alpha^\shortmid(t)&=\rho_0-\mathcal{K}\Theta\beta(t)-\frac{1}{2}\beta^\top
a \beta
\end{align}
with \emph{a} $ \in \mathbb{R}$ and $b \in \mathbb{R}^d$.\\

\noindent {\bf Proof of Proposition \ref{sol_0}}\\

Imposing that the processes $(X_t)_{t \geq 0}$ and $(\lambda_t)_{t
\geq 0}$ are affine as (\ref{SDE2}) and (\ref{inten}) so:
\begin{equation}
0=-(\rho_0+\rho_1\cdot
X)F(X,t)+F_t(X,t)+F_X(X,t)(\mathcal{K}(\Theta-X_t)) + \frac{1}{2}
\displaystyle\sum_{i,j}\frac{\partial F(X,t)}{\partial X_i
\partial X_j}(a_{ij}+b_{ij}\cdot X)
\end{equation}

Inserting $F(x,t)=e^{\alpha(t)+\beta(t)\cdot x}$ into the PDE
above and grouping the terms in $x$:
$$u(\cdot)x+v(\cdot)=0$$
Where
\begin{align}\label{ricati2}
u(\cdot)&= -\beta^\shortmid(t)+\rho_1-\mathcal{K}^\top
\beta(t)-\frac{1}{2}\beta(t)^\top
b \beta(t)\\
v(\cdot)&=\alpha^\shortmid(t)+\rho_0-\mathcal{K}\Theta\beta(t)-\frac{1}{2}\beta(t)^\top
a \beta(t)
\end{align}
Use the separation of variable technique to obtain that $\alpha$
and $\beta$ satisfy a Riccati equation with boundary condition
$\alpha(0)=0$ e $\beta(0)=w$.
\begin{flushright}
$\blacksquare$
\end{flushright}
%

%==================================

\section{Distribution of Doubly Stochastic Poisson Process}

For a non-homogeneous Poisson processes $(N_t)_{t \geq 0}$ with
intensity $(\lambda_t)_{t \geq 0}$, the probability of \emph{k}
occurrences within the interval $[t,T]$ is given by:
\begin{equation}\label{pdf}
P(N_T-N_t=k)= \frac{1}{k!}\mathbb{E}\left[\left(\int_t^T\lambda_u
du\right)^k \exp\left(-\int_t^T\lambda_u du\right)\right]
\end{equation}
$\mbox{for} \hspace{.3cm}  k=0,1,2,...$
\\ \\
Thus it is possible to state the important result:

\begin{theorem}\label{pdf2}
The Probability Distribution Function for a Doubly Stochastic
Poisson process, $(N_t)_{t \geq 0}$, within at interval $[t,T]$
when the intensity is an affine diffusion as (\ref{afim}) is given
by:
\begin{equation}
P(N_T-N_t=k)= \frac{1}{k!}\mathbb{E}\bigg[(\Lambda_t)^k e^{-
\Lambda_t}\bigg]=\frac{1}{k!}G_{\Lambda_t}^k(1)
\end{equation}

$\mbox{where} \hspace{.2cm} G_{\Lambda_t}^k(\mu)$ is the
\textit{k}-th derivative of Moment Generating Function for the
Hazard Process, $\Lambda_t$

\end{theorem}

\noindent {\bf Proof of Theorem \ref{pdf2}} \\

Imposing $f(X_t)\equiv 1$ and $V(X_t)\equiv \lambda_t$ on the left
hand side of equation (\ref{edp01}) it may be seen as the Laplace
Transform (or Moment Generating Function) for the Hazard Process,
$\Lambda_t:= \int_0^t \lambda_u du$.

Thus, using that if $G_X(\mu)$ is the MGF of $X$, then
\begin{equation}
 \frac{d^k G_X(\mu)}{d\mu^k} := G_{X}^k(\mu) = \mathbb{E}\left(X^k e^{\mu X}\right) \hspace{.1cm},
\end{equation}
and so the result follows.
\begin{flushright}
$\blacksquare$
\end{flushright}

In a particular case when the intensity follows an one-dimensional
Feller process we have:

\begin{theorem}\label{pdf3a}
The Probability Density Function (PDF) for a non-homogeneous
Poisson processes $(N_t)_{t \geq 0}$ within the interval $[t,T]$
when the intensity takes the form
\begin{equation}
d\lambda_{t} =\kappa  (\theta  -\lambda_{t}
)dt+\sigma\sqrt{\lambda_{t} } dW_{t} \end{equation}

is expressed by:

\begin{equation}\label{pdf4}
P(N_T-N_t=k)= \frac{1}{k!}\mathbb{E}\bigg[(\Lambda_t)^k e^{-
\Lambda_t}\bigg]=\frac{1}{k!}G_{\Lambda_t}^k(1)
\end{equation}

\begin{equation}
G_{\Lambda_t}(1)=\mathbb{E}\left(  e^{-\int_t^T \lambda_u
du}\right)= e^{\alpha(t,T)-\beta(t,T)\lambda_t}
\end{equation}

\noindent and

\begin{equation}\label{alpha}
\alpha(t,T)= \frac{2\kappa \theta }{\sigma ^{2} } \ln \left(
\frac{2\gamma \left( e^{(\gamma  +\kappa )/2} \right) }{\left(
\gamma  +\kappa \right) \left( e^{-\gamma (T-t)} -1\right)
+2\gamma } \right)
\end{equation}

%(5.22)\tab \tab \tab \tab
\begin{equation}\label{beta}
\beta(t,T)=\left[ \frac{2\mu\left(e^{-\gamma (T-t) }-1 \right)
}{\left( \gamma +\kappa \right) \left( e^{-\gamma (T-t)} -1\right)
+2\gamma } \right]
\end{equation}

\begin{equation}
\mbox{where} \hspace{.5cm} \gamma =\sqrt{\kappa^2 +2\sigma ^{2}}
\end{equation}
\end{theorem}
\noindent {\bf Proof of Proposition \ref{pdf3a}} \\
See appendix A

\begin{flushright}
$\blacksquare$
\end{flushright}

Since DSPPs processes are essentially Poisson processes, each
result for Poisson processes generally has a counterpart for DSPPs
processes. The following are some basic properties of a DSPP
process $(N_t)_{t \geq 0}$ driven by $(\lambda_t)_{t \geq 0}$. See
\citet{Bremaud} or \citet{Jones} for further properties.

\begin{proposition}\label{mom} The first two moments of $(N_t)_{t \geq 0}$
where the intensity is an one-dimensional Feller process are:
\end{proposition}
\hspace{-.2cm}    \textbf{1.}  $\mathbb{E}(N_t)$\hspace{.1cm}
\begin{equation}= \theta t+\frac{1-e^{-\kappa
    t}}{\kappa}(\lambda_0-\theta)\end{equation}\\

\hspace{-.70cm} \textbf{2}.\emph{Var}($N_t$)\hspace{.1cm} =
\begin{equation} \frac{2\theta t}{\kappa}\Bigl[(e^{-\kappa
t}+1)(\lambda_0-\theta)-2(\theta +
    \lambda_0)\Bigr]+
                   \frac{\sigma^2}{\kappa^3}\left(\frac{\theta e^{-\kappa
                   t}}{2} +\frac{4
                   e^{-\kappa}-5}{2}-\lambda_0 e^{-2\kappa t}\right) + \frac{\sigma^2
                   t}{\kappa}\left(\frac{3\theta-2\lambda_0}{\kappa}\right)
\end{equation}
\noindent {\bf Proof of Proposition \ref{mom}} \\

The first moment is obtained using the Laplace transform of
$\Lambda_t$:
$-\frac{\partial}{\partial\mu}\mathbb{E}\left(e^{-\mu\int_0^t
\lambda_s ds}\right)\Big|_{\mu=0}$
\newline

For the second moment we need an additional result:
\begin{equation}
\mbox{Var}\left(\Lambda_t\right)=
\frac{\partial^2}{\partial\mu^2}\mathbb{E}\left(e^{-\mu\int_0^t
\lambda_s ds}\right)\Big|_{\mu=0} -
\left(\frac{\partial}{\partial\mu}\mathbb{E}\left(e^{-\mu\int_0^t
\lambda_s ds}\right)\Big|_{\mu=0}\right)^2
\end{equation}
Substituting in $\mbox{Var}(N_t) =\mathbb{E}(\Lambda_t) +
\mbox{Var}(\Lambda_t)$ we obtain the result.

\begin{flushright}
$\blacksquare$
\end{flushright}

The stochastic nature of the intensity causes the variance of the
process to be greater than the variance of a homogeneous Poisson
process with the same expected intensity measure. This feature of
the DSPPs processes is referred, in the literature on point
processes and survival models, as \emph{overdispersion}.

%==============================
\section{Stationary distribution for $(N_t)_{t \geq 0}$}

Stationarity is a very important concept in time series analysis.
A stationarity assumption will allow us to estimate parameters
from point processes and make predictions. The characterization of
stationarity for a Point process relying on whether its intensity
process is stationary.\\

\begin{Definition} A point process $(N_t)_{t \geq 0}$ is stationary (or isotropic) if for all $A_1, \ldots, A_n \in \mathcal{F}$ e $h \in \mathbb{R}$
\begin{equation*}
(N_{A_1+h},\ldots,N_{A_n+h})\stackrel{d}{=}(N_{A_1},\ldots,N_{A_n})
\end{equation*}
\end{Definition}
Thus, the increments of $N_t$ are translation invariant
distribution with respect to any translation  $h \in
\mathbb{R}$.\\

The definition above can be written in a short form:
\begin{equation*}
\theta_h N\stackrel{d}{=}N,\hspace{.5cm} h \in \mathbb{R}
\end{equation*}
Where:
\begin{equation*}
\theta_hN_A := N_{A+h}
\end{equation*}
and $\stackrel{d}{=}$ means equally in distribution.  \\

%Na sequência apresentamos a seguinte proposição para a
%caracterização da estacionaridade do processo $N_t$:

\begin{proposition}\label{esta}
The Doubly Stochastic Poisson Process $(N_t)_{t \geq 0}$ is
stationary if its intensity $(\lambda_t)_{t \geq 0}$ is
stationary.
\end{proposition}
\noindent {\bf Proof of Proposition \ref{esta}} \\
%A transformada de Laplace de um processo pontual $N_t$ é definida
%por
%\begin{equation}\label{LN}
%L_N(f):= \mathbb{E}\left\{\exp\left [-\int_{\Omega}
%f(t)N(dt)\right] \right\}
%\end{equation}

The Laplace transform for a point process, $(N_t)_{t \geq 0}$,
with intensity process $(\lambda_t)_{t \geq 0}$ is given by:

\begin{equation}\label{LN1}
L_N(f)=\mathbb{E}\left\{\exp\left [-\int_{\Omega}
(1-e^{-f(t)})\lambda(t)(dt)\right] \right\}
\end{equation}
Thus, taking a $h \in \mathbb{R}$ with (\ref{LN1}):
\begin{align}
L_{\theta_hN} & =\mathbb{E}\left\{\exp\left [-\int_{\Omega} f(t)N(dt+h)\right] \right\} \nonumber \\
              & =\mathbb{E}\left\{\exp\left [-\int_{\Omega} f(t-h)N(dt)\right] \right\} \nonumber \\
              & =\mathbb{E}\left\{\exp\left [-\int_{\Omega} (1-e^{-f(t-h)})\lambda(t)(dt)\right] \right\} \nonumber \\
              & =\mathbb{E}\left\{\exp\left [-\int_{\Omega} (1-e^{-f(t)})\theta_h\lambda(t)(dt)\right] \right\}
              \label{LN2}
\end{align}
Therefore, if $(\lambda_t)_{t \geq 0}$ is stationary such that
$\theta_h \lambda\stackrel{d}{=}\lambda$, so together with
(\ref{LN2}) we obtain:
\begin{equation}
L_{\theta_hN}(f)=L_{N}(f)
\end{equation}
\begin{flushright}
$\blacksquare$
\end{flushright}

In a recent paper \citet{Glasserman} analyze the tail behavior,
the range of finite exponential moments, and the convergence to
stationarity in affine models, focusing on the class of canonical
models defined by \citet{DS}. According to \citet{Glasserman} the
one-dimensional Feller process has a stationary distribution so
the next step is determinate it.
%\begin{itemize}
%    \item Ornstein-Ulhenbeck (Vasicek)
%$$dX_t=\kappa(\theta-X_t)dt+\sigma dW_t$$
%    \item Feller (or Square Root, Cox-Ingersoll-Ross)
%$$dX_t=\kappa(\theta-X_t)dt+\sigma\sqrt{X_t} dW_t$$
%    \item Geometric Brownian Motion
%$$X_t:=ln S_t \hspace{1cm} dX_t=\left(\mu-\frac{1}{2}\sigma^2\right)dt+\sigma
%dW_t$$
%    \item Heston (Stochastic Volatility)
%    $$X_t:=ln S_t \hspace{1cm} dX_t=\left(\mu-\frac{1}{2}\sigma^2_t\right)dt+\sigma_t
%    dW^1_t$$ \\
%$$d\sigma_t=\kappa(\theta-\sigma_t)dt+\nu\sqrt{\sigma_t} dW^2_t$$\\
%$$<dW^1_t,dW^2_t>=\rho$$
%\end{itemize}

%Adicionalmente, deseja-se determinar de maneira analítica o
%comportamento de $N_t$ quando $t\longrightarrow \infty$. Desta
%forma, pode-se introduzir a seguinte definição:

\begin{proposition}\label{gama}
The stationary distribution of one-dimensional Feller process is
the Gamma distribution with parameters
$\alpha=\frac{2\kappa\theta}{\sigma^2}$ and
$\beta=\frac{2\kappa}{\sigma^2}$
\end{proposition}
\textbf{Proof of proposition \ref{gama}}\\
See Appendix B

\begin{flushright}
$\blacksquare$
\end{flushright}

\begin{Definition}\label{vago} Suppose that $\mu,\mu_1,\ldots$ are locally
finite measures defined on $\left( \Omega ,\mathcal{G},\mathbb{G},
\mathbb{P}\right)$, a necessary condition for \textbf{Vague
Convergence } of $\mu_n$ to $\mu$, in short $\mu_n
\stackrel{v}{\rightarrow}\mu$, is
$$ \mu_n f \rightarrow \mu f, \hspace{.1cm} \mbox{when} \hspace{.1cm} n\rightarrow
\infty,\hspace{.1cm} \mbox{for} \hspace{.1cm} f \in C_K^+(\Omega)
$$
\end{Definition}

Where: $C_K^+(\Omega)$ is a continuous function with compact
support  $f:\Omega\rightarrow \mathbb{R}_+$\\

From the definition \ref{vago} we shall present without proof the
well know result for convergence of point process:

\begin{result}\label{vago2} Let $N_1,N_2,\ldots$ defined on $\left( \Omega ,\mathcal{G},\mathbb{G},
\mathbb{P}\right)$ be a sequence of Point processes the following
results are equivalents:
\begin{enumerate}
    \item $N_t \stackrel{v}{\rightarrow}N_\infty$
    \item $N_tf \stackrel{v}{\rightarrow}N_\infty f,\hspace{.1cm} \mbox{for all} \hspace{.1cm} f \in
C_K^+(\Omega)$
    \item $\mathbb{E}\left(e^{N_tf}\right) \stackrel{v}{\rightarrow}\mathbb{E}\left(e^{N_\infty f}\right),\hspace{.1cm} \mbox{for all} \hspace{.1cm} f \in
C_K^+(\Omega)$
\end{enumerate}
\end{result}
\noindent {\bf Proof of Theorem \ref{vago2}} \\
\citet{kalenb}.
\begin{flushright}
$\blacksquare$
\end{flushright}

%Fazendo uso do Teorema \ref{vago2}, pode-se enunciar o seguinte
%resultado:

\begin{theorem}\label{vago3} For all $t\geq 1$, suppose that $(N_t)_{t \geq 0}$ is Poisson process
defined on $\left( \Omega ,\mathcal{G},\mathbb{G},
\mathbb{P}\right)$ with intensity $(\lambda_t)_{t \geq 0}$. If
$\lambda_t \stackrel{v}{\rightarrow}\lambda$ and $\lambda$ is
locally finite, thus $N_t \stackrel{v}{\rightarrow}N_\infty$,
where $N_\infty$ is a Poisson process with intensity $\lambda$
\end{theorem}

\noindent {\bf Proof of Theorem \ref{vago3}} \\
The Laplace transform for a point processes (\ref{LN1}) may be
written as $\mathbb{E}\left[e^{N_tf} \right]=e^{\lambda_t h}$, for
$f \in C_K^+(\Omega)$, where $h(t)=(1-e^{-f(t)})$. The Theorem
\ref{vago2} establish that $\lambda_t h \stackrel{v}{\rightarrow}
\lambda h$, so:
\begin{equation}
\mathbb{E}\left[e^{N_tf} \right]=e^{\lambda_t h}
\stackrel{v}{\rightarrow} e^{\lambda h}=
\mathbb{E}\left[e^{N_\infty f} \right]
\end{equation}
It follows from result \ref{vago2} that $N_t
\stackrel{v}{\rightarrow} N_\infty$.
\begin{flushright}
$\blacksquare$
\end{flushright}

From the above results it is possible to determine in a closed
form the stationary distribution of $(N_t)_{t \geq 0}$. We now
turn to one application of this result.

\begin{proposition} \label{Ninf}
The stationary distribution of $(N_t)_{t \geq 0}$ with intensity
$(\lambda_t)_{t \geq 0}$ given by $d\lambda_{t} =\kappa  (\theta
-\lambda_{t} )dt+\sigma\sqrt{\lambda_{t} } dW_{t}$ is the negative
Binomial distribution.
\end{proposition}
\noindent {\bf Proof of Proposition \ref{Ninf}}\\

The proof is straightforward when Proposition \ref{gama} is
considered together with Result \ref{vago2}

\begin{align}
\mathbb{P}(N_t=k) & = \int_0^\infty \mathbb{P}(N_t=k|\Lambda=\lambda)\Psi(\lambda)d\lambda \nonumber \\
              & =  \int_0^\infty \frac{(\lambda t)^k e^{-\lambda t}}{k!} \frac{\beta^\alpha}{\Gamma(\alpha)}e^{-\beta\lambda}\lambda^{\alpha-1}d\lambda \nonumber \\
              & =
              \frac{t^k}{k!}\frac{\beta^\alpha}{\Gamma(\alpha)}\frac{\Gamma(\alpha+k)}{(\beta+t)^{\alpha+k}}\int_0^\infty\frac{(\beta+t)^{\alpha+k}}{\Gamma(\alpha+k)}e^{-\lambda(\beta+t)}\lambda^{k+\alpha-1}d\lambda\nonumber
\end{align}
\begin{align}
\mathbb{P}(N_t=k)
              & =
              \frac{(\alpha+k-1)!}{(\alpha-1)!\alpha!}\left(\frac{t}{\beta+t}\right)^k\left(\frac{\beta}{\beta+t}\right)^\alpha \nonumber \\
              & \sim \mathcal{N}eg(\alpha,p) \hspace{.4cm} \mbox{with}
              \hspace{.5cm} p= \frac{\beta}{\beta+t}
              \nonumber
\end{align}

\begin{flushright}
$\blacksquare$
\end{flushright}

We have determined the stationary distribution for $(N_t)_{t \geq
0}$ when the intensity is given by an one-dimensional Feller
process. Closely related to this result is the determination of
the rate of convergence of $(N_t)_{t \geq 0}$ to its stationary
distribution $N_t^{\pi}$.

\begin{theorem}\label{speed} The Cox Process $(N_t)_{t
\geq 0}$ converge to its stationary distribution, $N_t^{\pi}$,
exponentially fast at rate $2 \kappa$.
\end{theorem}

\noindent {\bf Proof of Theorem \ref{speed}} It is known that the
transition density for the process $(\lambda_t)_{t \geq 0}$ when
it follows a Feller process is a non-central
$\chi^2$-distribution. From proposition \ref{gama} we have shown
that the process $(\lambda_t)_{t \geq 0}$ converge to the Gamma
distribution. According to \citet{Karlin}, it is possible rewrite
both distribution by its spectral representation:

\begin{equation}\label{spectral}
    p(t,x,y)= (2\kappa)^{2\kappa\theta}y^{2\kappa\theta-1}e^{2\kappa\theta y}\displaystyle\sum_{n=1}^{\infty}e^{-2\kappa
    t}L_n^{2\kappa-1}\left(\frac{\kappa x}{\sigma^2/2}\right)
      L_n^{2\kappa-1}\left(\frac{\kappa y}{\sigma^2/2}\right)\frac{\Gamma(n+1)}{\Gamma(n+\kappa\theta)}
\end{equation}

where $L_n^{2\kappa-1}(\cdot)$ is the Laguerre polynomial with
parameter $2\kappa-1$, and with the following property $L_0^{2\kappa-1}(\cdot)=1$. $\Gamma$ is the Gamma Function.\\

To describe the stationary distribution $p(y)$ by its spectral
representation we just need to set $n=0$ at equation (\ref{spectral}).\\

We decided to evaluate the convergence speed of $(\lambda_t)_{t
\geq 0}$ by looking at the convergence of its density functions
written in the spectral form:
\begin{equation}\label{speed2}
    |p(t,x,y)-p(y)|=\mathcal{O}(e^{-2\kappa t}),\forall x >0
    \hspace{.3cm} \mbox{when} \hspace{.3cm} t \rightarrow \infty
\end{equation}

Because according to the spectral representation, we have that for
$\forall n$, the only term inside the summation involving
\textit{t} is $e^{-2\kappa t}$.

Finally, the unconditional version of Cox process gives us:
\begin{equation}
P(N_T-N_t=k)=\int \left[\frac{1}{k!}\left(\Lambda_{t,T}\right)^k
\exp\left(-\Lambda_{t,T}\right)p(t,x,y)|\mathcal{F}_t\right]d\mathbb{P}
\end{equation}

Therefore,

\begin{equation}\label{speed3}
    |N_t-N_t^{\pi}|=\mathcal{O}(e^{-2\kappa t})    \hspace{.3cm} \mbox{when} \hspace{.3cm} t \rightarrow \infty
\end{equation}
\begin{flushright}
$\blacksquare$
\end{flushright}

\section{Model Application}

In section 2 we developed the requisite theory used to construct
the Cox Process as an affine function of the underlying state
variables. In every case, this relationship was subject to a given
parameter set. Unfortunately, the theory does not tell us anything
about the appropriate values that must be specified for this
parameter set. We must, therefore, turn to the econometric
literature to handle this important issue. Since the seminal paper
by Engle and Russell (1998) the modelling of financial point
process is an ongoing topic in the area of financial econometrics.
The financial point processes are associated with the random
arrival of specific financial trading events, such as
transactions, quote updates, limit orders or price changes
observable based on financial high-frequency data.

Moreover, it has been realized that the timing of trading events,
such as the arrival of particular orders and trades, and the
frequency in which the latter occur have information value for the
state of the market and play an important role in market
microstructure analysis.\\

Although the literature on the parametric estimation of point
processes (financial Point Process as well) is as large as the
theoretical literature, there is as yet no consensus as to the
best approach. Based on that, we propose a new technique to the
estimation of Cox Processes parameters. The methodology we will be
using, based on the Kalman filter, exploits the theoretical affine
relationship between the Cox Process and the state variables to
subsequently estimate the parameter set. The strength of this
approach is that Kalman filter is an algorithm that acts to
identify the underlying, and unobserved, state variables that
govern the Cox Process dynamics. Once the unobserved component is
filtered the Quasi-Maximum likelihood estimator will be able to
determine the model parameters.\\

In order to estimate parameters and to extract the unobservable
state variables we restrict the equation (\ref{pdf4}) to deal with
the probability of no arrivals, $P(N_T-N_t=0)$, within the
interval $(t-T)$. Therefore, the probability of no arrivals for
the Cox Process with Feller diffusion is

\begin{equation}\label{P0}
    P(N_T-N_t=0) = e^{\alpha(t,T)-\beta(t,T)\lambda_s}
\end{equation}
where $\alpha(t,T)$ and $\beta(t,T)$ were defined, respectively, in equations (\ref{alpha}) and (\ref{beta}).\\

Additionally it is possible linearize (\ref{P0}):

\begin{equation}\label{P1}
    \ln P(N_{t,T}=0) = \ln \alpha(t,T)-\beta(t,T)\lambda_s
\end{equation}

Thus, the measure equation is log-linear in $\lambda_t$ and it can
be written as:

\begin{equation}\label{Eq_med}
    \ln P(N_{t,T}=0) = \ln \alpha(t,T)-\beta(t,T)\lambda_s + \chi_s
\end{equation}
where: $$\chi_s \sim N(0, \textbf{R}_s)$$

The inclusion of an error term in equation (\ref{Eq_med}) is
motivated by the fact that the underlying intensity process may be
inadequate. If the true factor process is not a Feller process
equation (\ref{P0}) will be functionally misspecified and
estimates of $\lambda_t$ will be inferior. In this case the
probability of no arrivals within the interval $(t,T)$ implied by
the Feller process will systematically deviate from observed
arrivals. Therefore, in a correctly specified model the errors
$\chi_s$ should be serially and cross-sectionally uncorrelated
with mean zero.\\

It is known that the exact transition density
$p(t,x,y)=\mathbb{P}(X_t \in dy| X_{t-1}=dx)$ for the Feller
Process is the product of \emph{K} non-central $\chi^2$-densities.
Estimation of the unobservable state variables with an approximate
Kalman filter in combination with quasi-maximum-likelihood (QML)
estimation of the model parameters can be carried out by
substituting the exact transition density by a normal density:

$$ \lambda_s|\lambda_{s-1}\sim N(\mu_s,Q_s)$$
where $\mu_s$ and $Q_s$ are defined in such a way that the first
two moments of the approximate normal and the exact transition
density are equal. The moments are time varying and defined as:

\begin{equation}\label{Mom_1}
    \mu_s=\theta[1-\exp(-\kappa)]+\exp(-\kappa)\lambda_{s-1}
\end{equation}

and $Q_s$ is diagonal matrix with elements:

\begin{equation}\label{Mom_2}
    Q_s=\sigma^2\frac{1-\exp(-\kappa)}{\theta}\left(\frac{\theta}{2}[1-\exp(-\kappa)]+\exp(-\kappa)\lambda_{s-1}\right)
\end{equation}

Note that once we have represented the Cox Process with Feller
intensity into a State Space form we can, according to
\citet{Harvey}, to calculate the quasi log-likelihood:

\begin{equation}\label{QME}
    \log f(\textbf{y}_t|\textbf{x}_t;\pmb{\theta})=-\frac{1}{2}\log 2\pi(T-K)-\frac{1}{2}
    \displaystyle\sum_{t=K+1}^{T}\log|S_t|-\frac{1}{2}
    \displaystyle\sum_{t=K+1}^{T}\tilde{y}_t^{'}S_t^{-1}\tilde{y}_t
\end{equation}

where $T$ is the sample size; $K$ is the dimensional of the Feller
process.\\

%\vspace{.5cm}

\section{Empirical Analysis}

\subsection{Data description}
In this section, we apply our estimation technique to some
transaction data from the Brazilian Exchange
(BM\&FBOVESPA)\footnote{BM\&FBOVESPA is the fourth largest
exchange in the word in terms of market capitalization.
BM\&FBOVESPA has a vertically integrated business model with a
trade platform and clearing for equities, derivatives and cash
market for currency, government and private bonds. }.

The sample is formed by all submitted sell orders for BRL/USD FX
futures contract\footnote{Ticker: FUT DOLX08} traded in
BM\&FBOVESPA during October 2008 with expiry date of November 1,
2008. This FX contract is one of the most liquid FX contracts in
the emerging markets and the average volume of 300.000 traded
daily is significant even for developed markets. We have a total
of 535 records with trading occurring continuously from 10 am to 6
pm
exclusively through the electronic venue.\\

The BM\&FBOVESPA electronic trade system (GTS) uses the concept of
limit order book, matching orders by price/time priority. Lower
offer price take precedence over higher offers prices, and higher
bid price take precedence over lower bid prices. If there is more
than one bid or offer at same price, earlier bids and offers take
precedence over later bids and offers. The offers are recorded in
milliseconds allowing the highest precision for determination of
precedence criteria. We also observe that no two consecutive
orders arrive in the order book in an intervals smaller than 10
milliseconds, probably due to the internal network latency\\

While the probability of no arrivals are not themselves directly
observable, we can use the empirical frequency of orders arrivals
as proxy. To construct the empirical frequency we aggregate orders
sent within a same minute for every day during October and the
average value is taken. Thus, due to the 10 milliseconds network
latency, it is possible, at least theoretically, that 6.000 sells
orders could arrive in the order book in 60 seconds.

Therefore, the observable variable is constructed
as\footnote{where \# is the count measure within the interval
 $(t-s=60 \mbox{sec})$}:

\begin{equation}\label{var_obser}
    y_{(t-s=60\hspace{.1cm} \mbox{sec})}=\frac{\#(\mbox{Sell orders})}{6.000}
\end{equation}

\subsection{Estimation results}
In section we illustrate the properties of the Kalman filter for
estimating the parameters of Cox Process with Feller intensity.
Table 1 contains the parameter estimates
$\hat{\pmb{\theta}}=(\theta, \kappa, \sigma)$ . The estimated
standard deviations of errors - the square root of the diagonal
elements of H  - are also presented in Table 1. Standard errors
for the QML estimates are obtained as described in Hamilton (1994,
p.389).

\begin{table}[H]
 \begin{center}
\begin{tabular}{|c|c|c|c|c|}
  \hline
    & $\theta$ & $\kappa$ & $\sigma$ & Std. Error $\chi$ \\
  \hline
  Estimates & 0.065 & 0.0043 & 0.00267 & 0.0010 \\
  Std Error & 2.6E-05  & 1.7E-05 & 7.47E-08  & 2.68E-09 \\
  \hline
\end{tabular}
\caption{Estimated parameters of the Cox Process with affine
intensity applied for BRL/USD Futures sell orders data: (ticker
DOL FUT - maturity NOV09)}
 \end{center}
\end{table}

We obtain highly significant parameter estimates (at the 1\%
level). Additionally, the estimated standard deviation for the
measurement error $\chi$ is twice smaller than the diffusion
parameter $\sigma$.

%Figure 1 exhibits the model fitting in a in-sample analysis.

%
 \begin{figure}[!ht]
 \begin{center} %\vspace{-4.0cm}
\includegraphics[scale=0.5]{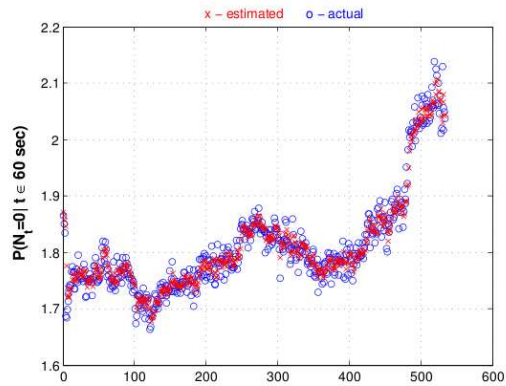}\vspace{-4.5cm}
\caption{Cox Process with Feller intensity fitted to BRL/USD
futures sell orders  (ticker DOL FUT - maturity NOV09)}
 \end{center}
\end{figure}

From figure 1 we can see the model flexibility in reacting to
different changes in the numbers of orders submitted. In this
case, the model seems to fit quite well to actual order flow and
the differences between actual and theoretical probability
indicate that, on average, the model tends to slightly
under-estimate actual occurrence of sell orders. In order to
assess the model fitted we conducted a diagnostic checking for
possible misspecifications based on the standardized residuals.
For a well specified model we need the residuals follows a white
noise process.
\begin{table}[H]
 \begin{center}
\begin{tabular}{|c|c|c|}
  \hline
  \emph{L} & Q(L) & p-value \\
  \hline
  5 & 8.36 & 0.13 \\
  10 & 13.47 & 0.19 \\
  15 & 17 & 0.31 \\
  \hline
\end{tabular}
 \end{center}
\begin{center}
\caption{Ljung-Box test for residual correlation}
\end{center}
\end{table}
 \vspace{.5cm}

Table 2 contains the Ljung-Box\footnote{The Ljung-Box statistic
tests the hypothesis that a process is serially uncorrelated.
Under the null hypothesis, the Q(L) statistic follows a Qui-Square
distribution with L degrees of freedom, where L is the maximum
number of temporal lags.} statistics up to the $15^{th}$ lag. The
LB test statistic allows us to not reject, at a high significance
level, the hypothesis that errors are white noise. In the absence
of dependence in errors the Feller process for intensity seems to
be a right choice to describe the observed sell orders behavior.

\section{Simulation results}

In this section, to analyze the performance of our estimation
algorithm we simulate various Cox Process outcomes using a known
parameter set and proceed to estimate the model parameters. This
simulation exercise is intended to indicate how effective this
technique is in terms of identifying parameters. The simulations
for the unobservable state variables have been drawn from the
noncentral Qui-Square distribution\footnote{The samples were drawn
from the Chi-Square distribution with $d$ degrees of freedom and
noncentrality parameter, $l$:

$$ c_s= \frac{\sigma^2(1-e^{-\kappa(t-s)})}{4\kappa}$$
$$d = \frac{4\kappa\theta}{\sigma^2}$$
$$l = \frac{\lambda_s e^{-\kappa(t-s)}}{c_s}$$}, $f_{\lambda}(d,l)$, and the measurement errors have been
simulated as normal random variables with zero means.\\

% This may not be an entirely sufficient number
%of iterations, but the procedure is rather lengthy and the results
%are primarily intended to give a general
%sense of the accuracy of the approach.\\

The simulation of a sample path for the Cox process with 500
elements, followed by application of the estimation algorithm, is
repeated 250 times. The following table summarize the results of
the simulation exercise for the Cox Process. The table reports the
true values, the mean estimate over the 250 simulations, the Mean
Quadratic Error, and the associated standard deviation of the
estimates.
\begin{table}[H]
\begin{center}
\begin{tabular}{|c|c|c|c|}
  \hline
   & $\theta$ & $\kappa$ & $\sigma$ \\
  \hline
  Actual Value & 0.04 & 0.2 & 0.05 \\
  Mean estimate & 0.04 & 0.3698 & 0.0448 \\
  MQE & 1.5303e-005 & 0.37 & 0.0023 \\
  Std. Error & 0.0039 & 0.52 & 0.0001 \\
  \hline
\end{tabular}
   \end{center}
\caption{A Analysis of the Kalman Filter estimates by  Monte Carlo
Simulation }
 \end{table}

 %\vspace{-8.0cm}
 \begin{figure}[!ht] %[b]
 \begin{center}
\includegraphics[scale=0.3]{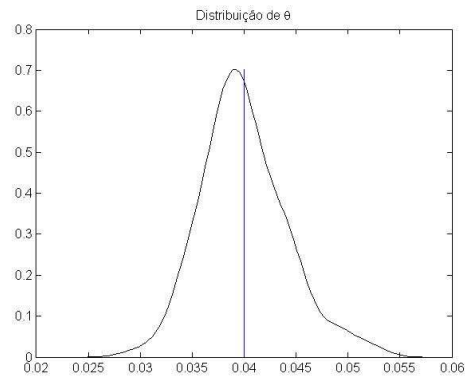}\vspace{-0.5cm} %{ofertas4.pdf}
 \end{center}
\caption{Empirical distribution of $\theta$}
\end{figure}

%\vspace{-15.5cm}
 \begin{figure}[!hb]
 \begin{center}
\includegraphics[scale=0.3]{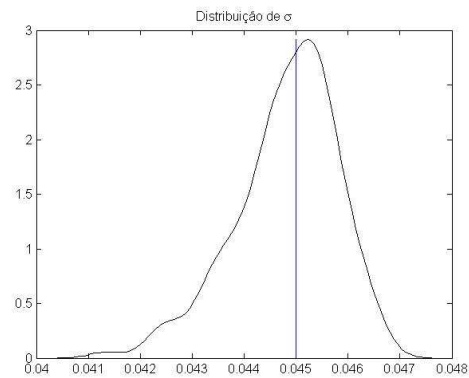}\vspace{-0.5cm} %{ofertas4.pdf}
\caption{Empirical distribution of $\sigma$}
 \end{center}
\end{figure}

 \begin{figure}[!ht] %[!ht]
 \begin{center}
\includegraphics[scale=0.3]{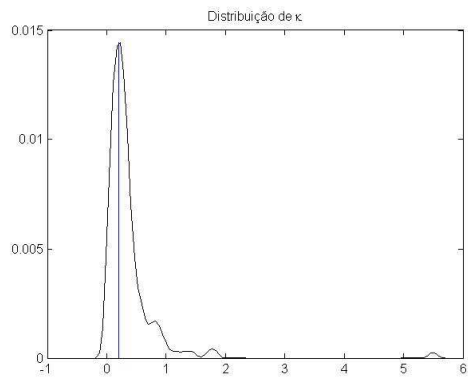}\vspace{-0.5cm} %{ofertas4.pdf}
\caption{Empirical distribution of $\kappa$}
 \end{center}
\end{figure}
%\vspace{-25.5cm}

%In reviewing Table 3, we can make the following specific
%observations.

%======================================
\clearpage
\section{Final Remarks}

Making use of the theoretical framework developed to model
interest rate term structure allow us to obtain the Probability
Density Function for a Double Stochastic Poisson Process (DSPP)
when the intensity process belong to a family of affine diffusion.
Furthermore, the stationary distribution for DSSP may be found
whenever the intensity process is also stationary. To illustrate
our results in this paper one of most common diffusion in interest
rate dynamics literature is assumed to drive the intensity
process, the one-dimensional Feller process. However the results
derived here are valid for any type of affine diffusion in
\emph{d}-dimension.\\

Finally this paper does not have an empirical focus, and these
results are primarily illustrative. The empirical analysis was
included to highlight that applied papers on estimating the
parameter set for the models examined in section 2 are
straightforward.

%======================================
\appendix
\section{ Proof of proposition \ref{pdf3}}
\noindent Specifying the form of $\mu(x,t)$ and $\sigma(x,t)$ into
the Theorem \ref{FK} we have:

\begin{equation} \label{edp}
\frac{\partial F}{\partial t}
+\kappa(\theta-\lambda)\frac{\partial F}{\partial \lambda}
+\frac{1}{2} \sigma^2\lambda \frac{\partial ^{2} F}{\partial
 \lambda^2} -\lambda\mu F=0
\end{equation}
With boundary condition $F(T,\lambda,\mu)=1$ \\ \newline

Using the results from proposition \ref{sol_0} together with
(\ref{edp})\footnote{Where subscript denotes partial derivative}:
\begin{align}
F_{\lambda}         & = -ABe^{-B \lambda} \nonumber \\
F_{\lambda\lambda}  & = AB^2e^{-B \lambda} \nonumber \\
F_{t}               & = A_te^{-B \lambda}- AB_t\lambda
e^{-B\lambda} \nonumber
\end{align}
Plugging the above results into (\ref{edp}):
\begin{equation}
\lambda\left(\frac{\sigma^2}{2}AB^2-AB_t+AB-A\right)=\kappa\theta
AB-A_t
\end{equation}
Since the left hand side is a function of $\lambda$, while the
right is independent of it, the following equations must be
satisfied:

\begin{equation}\label{ricati}
\left\{
\begin{array}{ccccccccc}
\frac{\sigma^2}{2}B^2 &-& B_t              &+& B &-& 1  &=& 0 \\
A_t                   &-& \kappa\theta AB  & &   & &    &=& 0 \\
\end{array}
\right.
\end{equation}

The first term of (\ref{ricati}) is Riccati equation with solution
$B(t,T)=v(t,T)/u(t,T)$ where $v(t,T)$ e $u(t,T)$ are solutions to
the following system: \footnote{$u$ e $v$ are functions of $t$ and
$T$, but $T$ is fixed, hence $v^{'}(t,T)$ denote the derivative
with respect to $t$}:

\begin{equation}\label{ricati2}
\left\{
\begin{array}{ccccccc}
\frac{\sigma^2}{2}v(t,T)&+&u^{'}(t,T)& &              &=&0 \\
u(t,T)                  &+&v^{'}(t,T)&-&\kappa v(t,T) &=&0 \\
\end{array}
\right.
\end{equation}
Let $\Delta=T-t$, so $\frac{\partial}{\partial t}=-\frac{d}{d
\Delta}$ and the system above may be written as:

\begin{equation}\label{ricati3}
\left\{
\begin{array}{ccccccc}
\frac{\sigma^2}{2}v(\Delta)&-&u^{'}(\Delta)& &              &=&0 \\
u(\Delta)                  &-&v^{'}(\Delta)&-&\kappa v(\Delta) &=&0 \\
\end{array}
\right.
\end{equation}
From the second term of (\ref{ricati3}) we have:

\begin{equation}\label{ricati4}
\left\{
\begin{array}{ccccccc}
u(\Delta)        &=&          &+&v^{'}(\Delta)&+&\kappa  v(\Delta)      \\
u^{'}(\Delta)    &=&          &+&v^{''}(\Delta)&+&\kappa v^{'}(\Delta)  \\
\end{array}
\right.
\end{equation}
Substituting into (\ref{ricati3}) and rewriting this in terms of
$\mathcal{D}$-operators:

\begin{equation}\label{ricati5}
 \left(\mathcal{D}^2 + \kappa\mathcal{D} -
\frac{\sigma^2}{2}\right)v(\Delta)=0
\end{equation}

Taking the roots of the quadratic equation (\ref{ricati5})
together with the boundary condition allow write the solution as:

\begin{equation}\label{ricati5}
v(\Delta)=e^{0.5(\gamma - \kappa)\Delta}-e^{0.5(-\gamma -
\kappa)\Delta}
\end{equation}

Substituting  into  (\ref{ricati4}) gives:
\begin{equation}\label{ricati6}
u(\Delta)=0.5(\gamma - \kappa) e^{0.5(\gamma -
\kappa)\Delta}-0.5(-\gamma + \kappa) e^{0.5(-\gamma -
\kappa)\Delta}
\end{equation}
Since $\Delta=T-t$, the solution of the Riccati equation is
obtained from (\ref{ricati5}) and (\ref{ricati6}):
\begin{align}
B(t,T)         & = v(\Delta)/u(\Delta) \nonumber \\
B(t,T)&=\left[ \frac{2\left(e^{-\gamma (T-t) }-1 \right) }{\left(
\gamma +\kappa \right) \left( e^{-\gamma (T-t)} -1\right) +2\gamma
} \right]\label{BtT}
\end{align}
Now consider equation (\ref{ricati}) with $T$ fixed, so $A(t,T)$
is a function of $t$ only. Hence:
\begin{align}
\frac{\partial A}{\partial t} &=\kappa\theta AB \nonumber \\
A(t,T) &=exp\left(-\kappa\theta\int_t^T B(s,T)ds\right) \nonumber
\end{align}
inserting $B(t,T)$ according to (\ref{BtT}) gives:

\begin{equation}
A(t,T)= \frac{2\kappa \theta }{\sigma ^{2} } \ln \left(
\frac{2\gamma \left( e^{(\gamma  +\kappa )/2} \right) }{\left(
\gamma  +\kappa \right) \left( e^{-\gamma (T-t)} -1\right)
+2\gamma } \right)
\end{equation}

\begin{flushright}
$\blacksquare$
\end{flushright}
%=======================================================
\section {Sketch of the proof of proposition \ref{gama}}

\noindent Let $p(s,t,x,y)$ be the transition density of
$\{\lambda_t: t \in \mathbb{R}_+\}$, for simplicity assume that
$\lambda_t=c \in \mathbb{R}_+ $ at instant $0$. From the Forward
Kolmogorov equation:
\begin{equation} \label{EFK}
\frac{\partial p}{\partial
t}(t,y)=\frac{1}{2}\frac{\partial^2}{\partial y^2}
\biggl[\sigma^2(y)p(t,y)\biggr] - \frac{\partial}{\partial}
\biggl[\mu(y)p(t,y)\biggr]
\end{equation}
The stationary distribution must satisfies $ \Psi(y)=\int
\Psi(x)p(t,x,y)dx \hspace{.3cm}\mbox{for all} \hspace{.3cm} t>0$.
When reached the stationary distribution is independent of $t$, so
$\frac{\partial p}{\partial t}(t,y)=0$. Combining it together with
(\ref{EFK}), we have\footnote{A rigorous proof might be found in
\citet{pinsky} p. 181 and 219}:
\begin{equation} \label{EFK2}
0=\frac{1}{2}\frac{\partial^2}{\partial
y^2}\biggl[\sigma^2(y)\Psi(y)\biggr] - \frac{\partial}{\partial}
\biggl[\mu(y)\Psi(y)\biggr]
\end{equation}
Let $s(y)=\exp \left(-\int_y \frac{2
\mu(\nu)}{\sigma^2(\nu)}d\nu\right)$ be the integrating factor, so
integrating (\ref{EFK2}) twice:
\begin{align}
\Psi(x) & = C_1\frac{S(x)}{s(x)\sigma^2(x)}+ C_2\frac{1}{s(x)\sigma^2(x)} \\
        & =m(x)[C_1S(x)+C_2]
\end{align}
Where: $S(x) = \int_x s(y)dy$\\

The constants $C_1$ e $C_2$ are determined such that $\Psi(x)$ be
a probability density function, so
$$\Psi(x)=\frac{1}{M}\frac{1}{\sigma^2(x)s(x)}dx \hspace{.4cm}
\forall \hspace{.2cm} x \in (l,r)$$

Thus, using the form of $s(x)$ and simplifying we have:
 $$\Psi(\lambda)=\beta^\alpha\frac{1}{\Gamma(\alpha)}e^{-\beta\lambda}\lambda^{\alpha-1}$$
In short, $\Psi(\lambda)=Gamma(\alpha,\beta)$

\bibliographystyle{plainnat_original}
\bibliography{Artigo_tese_v2}  % associado ao arquivo: 'bib_alan.bib'

%===================================
\end{document}